\documentclass{amsart}
\usepackage{amssymb,latexsym,epsfig,color}

\newcommand{\Z}{{\mathbb Z}^+ }
\newcommand{\R}{{\mathbb R}}

\newcommand{\set}{{\mathcal S}}
\newcommand{\M}{{\mathcal M}}

\newtheorem{theorem}{Theorem}[section]

\newtheorem{corollary}[theorem]{Corollary}
\newtheorem{lemma}[theorem]{Lemma}

\theoremstyle{definition}
\newtheorem{definition}[theorem]{Definition}
\newtheorem{remark}[theorem]{Remark}
\newtheorem{example}[theorem]{Example}
\newtheorem{algorithm}[theorem]{Algorithm}

\title [M-partitions]
{M-partitions: Optimal partitions of weight for one scale pan}

\author{Edwin O'Shea}
\thanks{The author was supported by a scholarship
under an agreement between the University of California and the
National University of Ireland, Cork during the bulk of this
work.}
\address{Department of Mathematics, University of
  Washington, Seattle, WA 98195-4350}
\email{oshea@math.washington.edu}
\date{\today}

\begin{document}

\begin{abstract}
 An $M$-partition of a positive integer $m$ is a partition with
 as few parts as possible such that any positive integer less
 than $m$ has a partition made up of parts taken from that
 partition of $m$.
 This is equivalent to partitioning a weight $m$ so as to be able
 to weigh any integer weight $l < m$ with as few weights as possible
 and only one scale pan.

 We show that the number of parts of an $M$-partition is a
 log-linear function of $m$ and the $M$-partitions of $m$
 correspond to lattice points in a polytope.
 We exhibit a recurrence relation for counting the number of
 $M$-partitions of $m$ and, for ``half'' of the positive
 integers, this recurrence relation will have a generating
 function. The generating function will be, in some sense,
 the same as the generating function for counting the number
 of distinct binary partitions for a given integer.
\end{abstract}

\maketitle
\section{Introduction}
\label{introduction}

Let $m$ be a positive integer and let
$\{\lambda_i \, : \, i \, =\, 0,1,\ldots,n \}$ be a finite
collection of, not necessarily distinct, positive integers with
$\lambda_0 \leq \lambda_1 \leq \cdots \leq \lambda_n$ and $m \, =
\, \lambda_0 + \lambda_1 + \cdots + \lambda_n$. In this case, we
say $m \, = \, \lambda_0 + \lambda_1 + \cdots + \lambda_n$ is a
{\em partition} of $m$ with $n+1$ {\em parts}. We will also refer
to the expression $\lambda_0 + \lambda_1 + \cdots + \lambda_n$ as
a partition. We call $\lambda_{i_0} + \lambda_{i_1} + \cdots +
\lambda_{i_k}$ a {\em subpartition} of the partition $m \, = \,
\lambda_0 + \lambda_1 + \cdots + \lambda_n$ if $\{\lambda_{i_0},\,
\lambda_{i_1}, \ldots, \lambda_{i_k} \}$ is a subcollection of
$\{\lambda_i \, : \, i \, =\, 0,1,\ldots,n \}$.

In~\cite{macmahon}, MacMahon called a partition $m \, = \,
\lambda_0 + \lambda_1 + \cdots + \lambda_n$ {\em perfect} if every
positive integer less than $m$ can be expressed uniquely as a
subpartition of $m \, = \, \lambda_0 + \lambda_1 + \cdots +
\lambda_n$. In this paper, we introduce partitions that are close
in spirit to MacMahon's. We maintain the subpartition property of
perfect partitions but drop the uniqueness constraint and we
demand that the number of parts in the partition be minimal.

\begin{definition} \label{Mdefn}
An {\bf M-partition} of $m$ is a partition $m=\lambda_0 +
\lambda_1 + \cdots + \lambda_n $ with $n$ being minimal such that
$\{\sum_{i \in I} {\lambda_i} \, : \, I \subseteq \{0,1, \ldots
,n\}\} = \{0,1,2, \ldots ,m\}$.
\end{definition}

We denote the set of all $M$-partitions for $m$ by $Mp(m)$. In
Section~\ref{parts} we will show that the number of parts in an
$M$-partition is a log-linear function of $m$ and that
$M$-partitions correspond to the lattice points in a certain
polytope. In particular, one can decide in polynomial time
whether a given partition is an $M$-partition or not.

\vspace{.2cm}

\noindent {\bf Theorem~\ref{numberofparts}.}
 {\em An $M$-partition of
 $m$ has precisely $\lfloor \log_2 m \rfloor + 1$ parts. }

\vspace{.2cm}

\noindent {\bf Theorem~\ref{inequalities}.}
 {\em 
 The partition
 $\lambda_0 + \lambda_1 + \cdots + \lambda_n$ is an $M$-partition
 if and only if $\lambda_i \leq  1 + \lambda_0 + \cdots + \lambda_{i-1}$
 for each $i \leq n$ and $2^n \leq
 \lambda_0 + \lambda_1 + \cdots + \lambda_n$. }

\vspace{.2cm}

In Section~\ref{algorithmsection} we develop algorithms for
generating $M$-partitions. These algorithms will be of great
benefit when proving the main result of
Section~\ref{countingsection} which is a recurrence relation for
counting the number of elements in $Mp(m)$, for each $m$. The
following is a special case of that recurrence relation.

\vspace{.2cm}

\noindent {\bf Theorem~\ref{simplecount}.}
 {\em 
 Let $m$ be a positive integer with $2^n + 2^{n-1} -1 \leq m < 2^{n+1}$
 for some positive integer $n$. Then }
 $
 | \, Mp(m) \, | = \sum_{i = \lfloor \frac{m}{2} \rfloor}
                       ^{2^n -1} { | \, Mp(i) \, | }.
 $

\vspace{.2cm}

In Section~\ref{simplifyingsection} we show that the recurrence
relation of Theorem~\ref{simplecount} is, in some sense, simultaneously
counting the number of $M$-partitions for an integer $m$ and counting
the number of distinct binary partitions for a given integer.

\vspace{.2cm}

\noindent {\bf Corollary~\ref{genfunction}.}
 {\em 
 If  $2^n + 2^{n-1} - 1 \leq m \leq 2^{n+1} - 1$ and
 $m = 2^{n+1} - 1 - k$
 then $| \, Mp(m) \, |$ equals the coefficient of
 $x^{ \lfloor \frac{k}{2} \rfloor }$
 in the generating function }
 $$
 (1-x)^{-1}\prod_{j=0}^{\infty}{(1-x^{2^{j}})^{-1}}.
 $$

In this paper $\Z$ will denote the positive integers and $m \in \Z$.
For every $r \in \R$ we denote by $\lceil r \rceil$
the smallest integer greater than or equal to $r$;
$\lfloor r \rfloor$ denotes largest integer less than or equal
to $r$. By $\log_2 m$ we mean the logarithm of $m$ base~$2$.

\section{The parts of an $M$-partition}
\label{parts}

We begin by investigating the subpartition property of
$M$-partitions. We define a weaker form of an $M$-partition by
dropping the minimality of parts constraint.
\begin{definition} \label{weakMdefn}
A {\bf weakM-partition} of $m$ is a partition $m=\lambda_0 +
\lambda_1 + \cdots + \lambda_n $ with $\{\sum_{i \in I}
{\lambda_i} \, : \, I \subseteq \{0,1, \ldots ,n\}\} = \{0,1,2,
\ldots ,m\}$.
\end{definition}

If $ m = \lambda_0 + \lambda_1 + \cdots + \lambda_n$ is a
$weakM$-partition of $m$ then we must have $\lambda_0 = 1$. If
$\lambda_1 \geq 3$ then it would not be possible to express $2$ as
a subpartition of $\lambda_0 + \lambda_1 + \cdots + \lambda_n$ and
so we must have $1 \leq \lambda_1 \leq 2$. In general, we have the
following bounds on the parts of a $weakM$-partition.
\begin{lemma} \label{biggun1}
  If $ m = \lambda_0 + \lambda_1 + \cdots + \lambda_n$ is a
  $weakM$-partition then
  $\lambda_i \leq 1 + \lambda_0 + \cdots + \lambda_{i-1}$
  for all $i \leq n$.
\end{lemma}
\begin{proof}
 Since $\lambda_i-1 < \lambda_i$ then $\lambda_i - 1$ can be
 expressed as $\lambda_i-1  = \sum_{j \in J} {\lambda_j}$ for some
 subset $J \subseteq \{0,1,\ldots,{i-1}\}$. Consequently,
 $\lambda_i - 1 \leq \lambda_0 + \lambda_1 + \cdots +
 \lambda_{i-1} $.
\end{proof}
\begin{lemma} \label{biggun2}
  If $ m = \lambda_0 + \lambda_1 + \cdots + \lambda_n$ is any
  partition with $\lambda_0 = 1$ and
  $\lambda_i \leq 1 + \lambda_0 \cdots + \lambda_{i-1} $
  for all $i \leq n$ then $\lambda_i \leq 2^i$ for all $i \leq n$.
\end{lemma}
\begin{proof}
 By assumption, $\lambda_0 = 1$. Proving by induction on $i$, assume
 $\lambda_k \leq 2^k$ for all $k \leq i-1$. We are given that
 $\lambda_i \leq 1 + \lambda_0 \cdots + \lambda_{i-1} $
 and so by the induction hypothesis, we have
 $\lambda_i \leq 1 + \sum_{k=0}^{i-1} {2^k} = 2^i$.
\end{proof}

The upshot of Lemma~\ref{biggun1} and Lemma~\ref{biggun2} is a
lower bound on the number of necessary parts in a
$weakM$-partition.
\begin{corollary} \label{necessaryparts}
  If $ m = \lambda_0 + \lambda_1 + \cdots + \lambda_n$ is a
  $weakM$-partition then $n \geq \lfloor \log_2 m \rfloor$.
\end{corollary}
\begin{proof}
 If $ m = \lambda_0 + \lambda_1 + \cdots + \lambda_n$
 is a $weakM$-partition of $m$ then
 $ m = \lambda_0 + \lambda_1 + \cdots + \lambda_n \leq 2^{n+1}-1 < 2^n$.
 This implies that $\lfloor \log_2 m \rfloor < n+1$. Since
 $\lfloor \log_2 m \rfloor$ is an integer then it is no more than $n$.
\end{proof}
\begin{remark}
 Lemma~\ref{biggun1} and Lemma~\ref{biggun2} apply equally to
 $M$-partitions since every $M$-partition is a $weakM$-partition.
 Corollary~\ref{necessaryparts} provides a lower bound for the
 minimality of parts criterion of $M$-partitions.
\end{remark}
It is well known that every postive integer has a unique binary
representation and this has the following implication for
$weakM$-partitions.
\begin{lemma} \label{binaryalgorithm}
  The partition $1 + 2 + 4 + \cdots + 2^n$ is a $weakM$-partition
  of $2^{n+1} - 1$.
\end{lemma}
\begin{remark} \label{complement}
 In order to show that a partition of $m$ is a $weakM$-partition it
 is sufficient to show that for all $l \leq \lceil \frac{m}{2} \rceil$
 there is some $J \subseteq \{0, 1 , \ldots , n \}$ with
 $\sum_{j \in J}{\lambda_j} = l$, since
 $m - l = \sum_{j \in J^c}{\lambda_j}$
 where $J^c$ is the complement of $J$.
\end{remark}
The following algorithm shows that the lower bound presented for
the number of parts in Corollary~\ref{necessaryparts} is
sufficient.
\begin{algorithm} \label{algorithm1}
  There exists a $weakM$-partition of $m$ with
  $\lfloor \log_2 m \rfloor + 1 $ parts.
\end{algorithm}
\begin{proof}
 Let $n =  \lfloor \log_2 m \rfloor$ and list the $n+1$ integers
 $2^0, 2^1, 2^2, \ldots, 2^ {n-1}, m-(2^n - 1)$ in
 increasing order and set a one-to-one correspondence with
 $\lambda_0, \lambda_1, \ldots, \lambda_n$. Then
 $m = \lambda_0 + \lambda_1 + \cdots + \lambda_n$ is a partition and
 we claim that every $l < m$ can be expressed a subpartition of
 this partition.

 If $m = 2^{n+1} - 1$ then by Lemma~\ref{binaryalgorithm}
 we are done. Otherwise, by Corollary~\ref{necessaryparts},
 $m \leq 2^{n+1} - 2$ and so $\lceil \frac{m}{2} \rceil \leq 2^n -1$.
 Since the parts of $2^n - 1 = 2^0 + 2^1+ 2^2+ \cdots + 2^{n-1}$ are
 all parts of the partition given then, combining
 Lemma~\ref{binaryalgorithm} with Remark~\ref{complement}, we see
 that $2^0, \, 2^1, \, 2^2, \ldots,  2^ {n-1}, \, m-(2^n - 1)$ are
 the parts of a $weakM$-partition of $m$.
\end{proof}
\begin{example} \label{mequals53}
 Let $m=53$. Using Algorithm~\ref{algorithm1} we have the
 $weakM$-partition $53 = 1+2+4+8+16+22$.
\end{example}
The first main result of this section is that the above algorithm
describes a way to find an $M$-partition for any $m$.
\begin{theorem} \label{numberofparts}
  An $M$-partition of $m$ has precisely
  $\lfloor \log_2 m \rfloor + 1$ parts.
\end{theorem}
\begin{proof}
 Corollary~\ref{necessaryparts} asserts that at least
 $\lfloor \log_2 m \rfloor + 1$ parts are needed for an $M$-partition
 of $m$. Algorithm~\ref{algorithm1} tells us that this is sufficient.
\end{proof}

\begin{example} \label{countingexample}
\begin{align*}
 Mp(7) & = \{1+2+4\},\\
 Mp(8) & = \{1+1+2+4, \, 1+1+3+3, \, 1+2+2+3\},\\
 Mp(9) & = \{1+1+2+5, \, 1+1+3+4, \, 1+2+2+4, \, 1+2+3+3\},\\
 Mp(10)& = \{1+1+3+5, \, 1+2+2+5, \, 1+2+3+4\},\\
 Mp(11)& = \{1+1+3+6, \, 1+2+2+6, \, 1+2+3+5, \, 1+2+4+4\},\\
 Mp(12)& = \{1+2+3+6, \, 1+2+4+5\} , \, Mp(13) = \{1+2+3+7, \, 1+2+4+6\},\\
 Mp(14)& = \{1+2+4+7\} , \, Mp(15) = \{1+2+4+8\}.
\end{align*}
 You will need $5$ parts for each $M$-partition of $16$ and there
 $12$ such $M$-partitions.
\end{example}

At first sight, it appears that deciding whether a partition is a
$weakM$-partition or not could be an arduous endeavor. However, we
have a relatively painless way of deciding so which avoids
checking that the subpartition property holds for every $l <m$.
\begin{lemma} \label{weakMinequalities}
  The partition $\lambda_0 + \lambda_1 + \cdots + \lambda_n$ is a
  $weakM$-partition if and only if
  $\lambda_i \leq  1 + \lambda_0 + \cdots + \lambda_{i-1}$ for each
  $i \leq n$.
\end{lemma}
\begin{proof}
 The ``only if'' follows from Lemma~\ref{biggun1}. Conversely let $\set_n$
 be the set of all partitions with $n+1$ parts that satisfy
 $\lambda_i \leq 1 + \lambda_0 + \cdots + \lambda_{i-1}$
 for each $i \leq n$. We will argue the ``if'' by showing
 that $\set_n$ is contained in the set of $weakM$-partitions with
 $n+1$ parts. We will do so by induction on $n$.

 It is clear that $\set_0 = \{1\}$. Assume the induction hypothesis
 on $\set_i$ for all $i \leq n-1$. Let $\lambda_0 + \lambda_1 + \cdots +
 \lambda_n$ be a partition in $\set_n$ and let
 $l < \lambda_0 + \lambda_1 + \cdots + \lambda_n$. We need to
 show that $l$ can be expressed as a subpartition of
 $\lambda_0 + \lambda_1 + \cdots + \lambda_n$. Note that
 $\lambda_0 + \lambda_1 + \cdots + \lambda_{n-1}$ is in
 $\set_{n-1}$ and so by our induction hypothesis if
 $l \leq \lambda_0 + \lambda_1 + \cdots + \lambda_{n-1}$
 then there is nothing to show. Hence, we only need concern
 ourselves with $\lambda_{n-1} < l < \lambda_n$ and
 $l > \lambda_n$.

 If $\lambda_{n-1} < l < \lambda_n$ then
 $l - \lambda_{n-1} <  \lambda_n - \lambda_{n-1}$. By virtue
 of $\lambda_0 + \lambda_1 + \cdots + \lambda_n$ being in $\set_n$ we have
 $\lambda_n - \lambda_{n-1} \leq (1 + \lambda_0 + \cdots + \lambda_{n-1})
                                 - (\lambda_{n-1})$
 and so
 $l - \lambda_{n-1} <  \lambda_0 + \lambda_1 + \cdots +
 \lambda_{n-2}$. But the partition $\lambda_0 + \lambda_1 +
 \cdots + \lambda_{n-2}$ is in $\set_{n-2}$ and so $l$ can be expressed
 in terms of a subpartition of $\lambda_0 + \lambda_1 +
 \cdots + \lambda_{n-1}$ which is a subpartition of
 $\lambda_0 + \lambda_1 + \cdots + \lambda_{n-1} + \lambda_n$.
 Similarly, since $l < \lambda_0 + \lambda_1 + \cdots + \lambda_n$,
 $l > \lambda_n$ implies that
 $0 < l - \lambda_n < \lambda_0 + \lambda_1 \cdots +
 \lambda_{n-1}$.
 By our inductive hypothesis, $l - \lambda_n$ can be expressed a
 subpartition of $\lambda_0 + \lambda_1 \cdots + \lambda_{n-1}$
 and so $l$ can be expressed as a subpartition of
 $\lambda_0 + \lambda_1 \cdots + \lambda_n$.
\end{proof}

The second main result of this section is that there is an
efficient way of deciding whether a given partition is an
$M$-partition or not. This is achieved by a polyhedral
characterization of $M$-partitions.
\begin{theorem} \label{inequalities}
  The partition $\lambda_0 + \lambda_1 + \cdots + \lambda_n$ is an
  $M$-partition if and only if
  $\lambda_i \leq 1 + \lambda_0 + \cdots + \lambda_{i-1}$
  for each $i \leq n$ and
  $2^n \leq \lambda_0 + \lambda_1 + \cdots + \lambda_n$.
\end{theorem}
\begin{proof}
 The ``only if'' follows from Lemma~\ref{biggun1} and
 Theorem~\ref{numberofparts}. As for the converse we need to show
 that
 $n = \lfloor \log_2(\lambda_0 + \lambda_1 + \cdots + \lambda_n)
 \rfloor$ and that $\lambda_0 + \lambda_1 + \cdots + \lambda_n$ is
 a $weakM$-partition.

 From Lemma~\ref{biggun2} we have that
 $\lambda_0 + \lambda_1 + \cdots + \lambda_n < 2^{n+1}$ and, by
 assumption, we have
 $2^n \leq \lambda_0 + \lambda_1 + \cdots + \lambda_n$.
 Therefore, the partition $\lambda_0 + \lambda_1 + \cdots + \lambda_n$
 has the desired number of parts. From Lemma~\ref{weakMinequalities} we have
 $\lambda_0 + \lambda_1 + \cdots + \lambda_n$ is a $weakM$-partition.
\end{proof}

An important consequence of Theorem~\ref{inequalities} is that
$M$-partitions are both built upon, and can be extended to, other
$M$-partitions.
\begin{corollary} \label{extendedM}
  Let $m = \lambda_0 + \lambda_1 + \cdots + \lambda_n$ be an
  $M$-partition. Then $\lambda_0 + \lambda_1 + \cdots
  + \lambda_j$ is an $M$-partition for all $j \leq n$.
  Also, if $r \in \Z$ then the partition
  $ m + r = \lambda_0 + \lambda_1 + \cdots + \lambda_n + r$
  is an $M$-partition of $m + r$ if and only if
  $\lambda_n \leq r$, $r \leq m + 1$ and $2^{n+1} \leq m + r$.
\end{corollary}
\begin{proof}
 Since
 $\lambda_i \leq 1 + \lambda_0 + \cdots + \lambda_{i-1}$
 for each $i \leq n$ then
 $\lambda_i \leq (1 + \lambda_0 + \cdots + \lambda_{i-2} )+ \lambda_{i-1}
            \, \leq \, (1 + \lambda_0 + \cdots + \lambda_{i-2} )
            +    (1 + \lambda_0 + \cdots + \lambda_{i-2} )
            \, = \, 2 (1 + \lambda_0 + \cdots + \lambda_{i-2} )
 $.
 Continuing in this fashion we can see that
 $\lambda_i \leq 2^{i-j-1} (1 + \lambda_0 + \cdots + \lambda_j )$
 for all $ i > j$. Since
 $ 2^n \leq  \lambda_0 + \lambda_1 + \cdots + \lambda_n $
 then
 $2^n \leq 2^{n-j}(\lambda_0 + \cdots + \lambda_j) +
 (2^{n-j} -1)$.
 Therefore, $2^j \leq \lambda_0 + \cdots + \lambda_j$ since
 $\frac{(2^{n-j} -1)}{2^{n-j}} < 1$.
 Since $\lambda_i \leq 1+ \lambda_0 + \cdots + \lambda_{i-1}$
 for each $i \leq j$, then
 $\lambda_0 + \lambda_1 + \cdots + \lambda_j$ is an $M$-partition
 for all $j \leq n$.

 Next, $ m + r = \lambda_0 + \lambda_1 + \cdots + \lambda_n + r$
 is a partition which, by definition, means $\lambda_n \leq r$.
 We assumed $m = \lambda_0 + \lambda_1 + \cdots + \lambda_n$ to be
 an $M$-partition so, by Theorem~\ref{inequalities}, both
 $r \leq m + 1$ and $2^{n+1} \leq m + r$ are necessary and sufficient
 for our claim.
\end{proof}
\begin{remark} \label{extendedMremark}
 An important reformulation of the extension statement in
 Corollary~\ref{extendedM} is the following: Let $m \in \Z$
 with $n = \lfloor \log_2 m \rfloor$ and let $m^{(1)} < m$. Then
 $m = \lambda_0 + \lambda_1 + \cdots + \lambda_{n-1} + ( m - m^{(1)})$
 is an $M$-partition if and only if
 $\lambda_{n-1} \leq m - m^{(1)}$,
 $m^{(1)} = \lambda_0 + \lambda_1 + \cdots + \lambda_{n-1}$
 is an $M$-partition and $ m - m^{(1)} \leq m^{(1)} + 1$.
\end{remark}

For the rest of this exposition, in light of
Theorem~\ref{numberofparts} and Theorem~\ref{inequalities}, all
partitions will be $M$-partitions unless otherwise stated, and $n$
will always refer implicitly to some $m$ via $n = n(m) := \lfloor
\log_2 m \rfloor$.


\section{Algorithms for generating $M$-partitions}
\label{algorithmsection}

In this brief section we give two more algorithms for generating
$M$-partitions. These algorithms, in addition to
Algorithm~\ref{algorithm1}, will assist us in attaining an exact
count for the number of $M$-partitions of $m$ for all $m \in \Z$.

\begin{algorithm} \label{algorithm2}
  Letting $m \in \Z$, assign $\lambda_n  =  \lceil \frac{m}{2}
  \rceil$ and recursively define
  $$\lambda_i  = \lceil \frac{m - (\lambda_n + \lambda_{n-1} +
     \cdots + \lambda_{i+1})}{2} \rceil
  $$
  for all non-negative $i < n$. Then
  $m = \lambda_0 + \lambda_1 + \cdots + \lambda_n$
  is a partition of $m$.
\end{algorithm}
\begin{proof}
 By construction,
 $\lambda_0 \leq \lambda_1 \leq \cdots \leq \lambda_n$. Let $T_n$
 be the statement
 ``if $m \in \Z$ with $n = \lfloor \log_2 m \rfloor$ then
 $m = \sum_{i=0}^n{\lambda_i}$. '' We will show by induction
 that $T_n$ is true for all $n$.

 The statement $T_0$ is true since $1 = \lceil \frac{1}{2}
 \rceil$. Assume that $T_{n-1}$ is true. Let $m \in \Z$ with
 $n = \lfloor \log_2 m \rfloor$. Then
 $\lambda_n  =  \lceil \frac{m}{2} \rceil$ and so
 $m - \lambda_n = \lfloor \frac{m}{2} \rfloor$. But
 $ \log_2 \lfloor \frac{m}{2} \rfloor = n-1$ and so
 $\lambda_0 + \lambda_1 + \cdots + \lambda_{n-1}$ is a partition
 of $\lfloor \frac{m}{2} \rfloor$ since $T_{n-1}$ is assumed to
 be true. Hence $\lambda_0 + \lambda_1 + \cdots + \lambda_{n-1} + \lambda_n
 = \lfloor \frac{m}{2} \rfloor + \lceil \frac{m}{2} \rceil = m$.
\end{proof}
\begin{corollary} \label{partition2}
  The partition $m = \lambda_0 + \lambda_1 + \cdots + \lambda_n$
  given by Algorithm~\ref{algorithm2} is an $M$-partition.
\end{corollary}
\begin{proof}
 Since $n = \lfloor \log_2 m \rfloor$ then
 $2^n \leq \lambda_0 + \lambda_1 + \cdots + \lambda_n$.
 Since $m = \lambda_0 + \lambda_1 + \cdots + \lambda_n$ then
 $\lambda_i = \lceil \frac{\lambda_0 + \lambda_1 +
  \cdots + \lambda_i} {2} \rceil
  <
  \frac{\lambda_0 + \lambda_1 + \cdots + \lambda_i} {2} + 1
 $. Therefore,
 $\lambda_i \leq  1 + \lambda_0 + \cdots + \lambda_{i-1}$ for all
 $i \leq n$. By Theorem~\ref{inequalities} the partition described
 is an $M$-partition.
\end{proof}

Algorithm~\ref{algorithm1} and Algorithm~\ref{algorithm2} provide
$M$-partitions with $n+1$ parts for all $m$ such that $2^n \leq m
< 2^{n+1}$. The next algorithm offers an $M$-partition for $m$ if
there is the further restriction that $2^n \leq m \leq
2^n+2^{n-1}-2$. The need for such a special case will become
apparent in Section~\ref{countingsection}.
\begin{algorithm} \label{algorithm3}
  Let $m \in \Z$ with $2^n \leq m \leq 2^n+2^{n-1}-2$. Define
  $\lambda_i = 2^i $ for all $i \leq n-2$,  $\lambda_{n-1} =
  \lfloor \frac{m - (2^{n-1} -1) }{2} \rfloor$ and $\lambda_n =
  \lceil \frac{m - ( 2^{n-1} -1) }{2} \rceil$.
  Then $m = \lambda_0 + \lambda_1 + \cdots + \lambda_n$ is an
  $M$-partition.
\end{algorithm}
\begin{proof}
 It is clear that this algorithm provides a partition of $m$.
 By Theorem~\ref{numberofparts} the partition has the desired
 number of parts. All we need show is that
 $m = \lambda_0 + \lambda_1 + \cdots + \lambda_n$ is a
 $weakM$-partition.

 By Remark~\ref{complement} all we need show that every
 $l \leq \lceil \frac {m}{2} \rceil \leq 2^{n-1}+2^{n-2}-1$
 can be expressed as a subpartition of
 $m = \lambda_0 + \lambda_1 + \cdots + \lambda_n$.
 If $l \leq 2^ {n-1}-1$ then Lemma~\ref{binaryalgorithm}
 applies and $l$ can be expressed as a subpartition of
 $\lambda_0 + \lambda_1 + \cdots + \lambda_{n-2}$.
 Alternatively, suppose $2^{n-1} \leq l \leq
 \lceil \frac{m}{2} \rceil \leq 2^{n-1}+2^{n- 2}-1$.
 By our restrictions on $m$ and our choice of $\lambda_n$
 we have $\lambda_n \geq 2^{n-2}$ and hence,
 $l - \lambda_n \leq 2^{n-1} - 1$. By
 Lemma~\ref{binaryalgorithm}, $l - \lambda_n$ can be expressed
 as a subpartition of
 $\lambda_0 + \lambda_1 + \cdots + \lambda_{n-2}$ and thus
 $l$ can be expressed as a subpartition of
 $m = \lambda_0 + \lambda_1 + \cdots + \lambda_n$.
\end{proof}

\vspace{.2cm} \noindent
 {\bf Example~\ref{mequals53} continued.}
 Algorithm~\ref{algorithm2} yields the $M$-partition
 $53 = 1+2+3+7+13+27$. Algorithm~\ref{algorithm3} produces
 the partition $53=1+2+4+8+9+9$ but this is not an $M$-partition
 as we have no way of expressing $16$ as a subpartition.

\section{Counting the number of elements in the set $Mp(m)$}
\label{countingsection}

For each $m \in \Z$ define $Mp(m)$ to be the set of all
$M$-partitions of $m$. By Corollary~\ref{extendedM} and
Remark~\ref{extendedMremark} we know that every $M$-partition must
be constructed upon another of one less part. Letting $a_m :=
|Mp(m)|$ we construct a recurrence relation for $a_m$ by way of
finding sharp bounds on the largest part of an $M$-partition of
$m$.

Fix $m \in \Z$. Let $m^{(1)} \in \Z$ be any integer whose
$M$-partitions can be extended to an $M$-partition of $m$ in the
sense of Remark~\ref{extendedMremark}. Similarly, for each such
$m^{(1)}$, let $m^{(12)} \in \Z$ be any integer whose
$M$-partitions can be extended to an $M$-partition of $m^{(1)}$.

\vspace{.2cm}

\noindent {\bf Remark~\ref{extendedMremark} continued.}
 The number of $M$-partitions of $m$, $a_m$ equals the cardinality of
 the set of partitions given by
 \[
  \left \{
   \, \lambda_0 + \lambda_1 + \cdots + \lambda_{n-1} \, :
    \begin{array}{ll}
    & m^{(1)} < m , \, \lambda_{n-1} \leq m - m^{(1)},
    \,  m - m^{(1)} \leq m^{(1)} + 1 \\
    & \textup{and}\, \lambda_0 + \lambda_1 + \cdots + \lambda_{n-1} \,
    \textup{is an $M$-partition of} \, m^{(1)}
    \end{array}
   \right \}.
 \]

\vspace{.2cm}

We now turn our attention to determining what values these
$m^{(1)}$'s can take on for a given $m$. We do so by determining
sharp bounds on the largest part of an $M$-partition of $m$.
\begin{lemma} \label{lowerbound1}
 Let $m = \lambda_0 + \lambda_1 + \cdots + \lambda_n$ be an
 $M$-partition.
 Then $$ \lceil \frac{m - 2^{n-i+1} + 1}{i} \rceil \leq \lambda_n.$$
\end{lemma}
\begin{proof}
 By Lemma~\ref{biggun2} we have $\lambda_i \leq 2^i$ for all $i
 \leq n$. Since $\lambda_{i-1} \leq \lambda_i$ for all $i \leq n-1$
 then
 $ i \lambda_n \geq \lambda_{n-i+1} + \cdots + \lambda_n
 = m - ( \lambda_0 + \lambda_1 + \cdots + \lambda_{n-i})
 \geq m - 2^{n-i+1} + 1$.
 Hence, $\lambda_n \geq \lceil \frac{m - 2^ {n-i+1}+ 1}{i} \rceil $
\end{proof}
\begin{remark} \label{onlytwolowerbounds}
 It is unnecessary to consider all of the bounds in
 Lemma~\ref{lowerbound1} -- we only need consider the bounds
 given by $i=1$ and $i=2$.
 When $2^n + 2^{n-1} -1 \leq m \leq 2^{n+1}-1$ then
 $m - 2^n +1 \geq \lceil \frac{m - 2^{n-i+1} + 1}{i} \rceil$
 for all $i \leq n$. If $2^n \leq m \leq 2^n + 2^{n-1} -2$
 then $ \lceil \frac{m - 2^{n-1} + 1}{2} \rceil  \geq
 \lceil \frac{m - 2^{n-i+1} + 1}{i} \rceil$ for all $i \leq n$.
\end{remark}
\begin{lemma} \label{bound1}
 Let $m = \lambda_0 + \lambda_1 + \cdots + \lambda_n$ be an
 $M$-partition. Then
 $$ max \{m - 2^n +1, \, \lceil \frac{m - 2^{n-1} + 1}{2} \rceil \}
 \, \leq \, \lambda_n \leq \, \lceil \frac{m}{2} \rceil.$$
 Furthermore, all three bounds are sharp.
\end{lemma}
\begin{proof}
 If $\lambda_n > \lceil \frac{m}{2} \rceil$ then
 $\lambda_0 + \lambda_1 + \cdots + \lambda_{n-1}
 < m - \lceil \frac{m}{2} \rceil = \lfloor \frac{m}{2} \rfloor
 \leq   \lceil \frac{m}{2} \rceil < \lambda_n $.
 This implies $\lfloor \frac{m}{2} \rfloor$ cannot be
 expressed as a subpartition of
 $m = \lambda_0 + \lambda_1 + \cdots + \lambda_n$ which
 contradicts $m = \lambda_0 + \lambda_1 + \cdots + \lambda_n$
 being an $M$-partition. Hence $\lambda_n \leq \lceil \frac{m}{2} \rceil$.
 The lower bounds follow from Lemma~\ref{lowerbound1} and
 Remark~\ref{onlytwolowerbounds}. Algorithm~\ref{algorithm1},
 Algorithm~\ref{algorithm2} and Algorithm~\ref{algorithm3} insure
 that all three bounds can be attained for any given $m$.
\end{proof}
\begin{corollary} \label{bound1flip}
 Let $m = \lambda_0 + \lambda_1 + \cdots + \lambda_n$ be an
 $M$-partition. Then
 $$
 \lfloor \frac{m}{2} \rfloor
 \leq \lambda_0 + \lambda_1 + \cdots + \lambda_{n-1}
 \leq min\{ \lfloor \frac{m + 2^{n-1} - 1}{2} \rfloor ,
 2^n - 1 \}.
 $$
\end{corollary}
For a given $m$, we can restate Corollary~\ref{bound1flip} in
terms of the $m^{(1)}$'s and in turn for the $m^{(12)}$'s of each
such $m^{(1)}$.
\begin{corollary} \label{bound1flipnew}
 Let $m \in \Z$. Then
 \begin{equation} \label{strata_1_bound}
 \lfloor \frac{m}{2} \rfloor
 \leq m^{(1)}
 \leq min\{ \lfloor \frac{m + 2^{n-1} - 1}{2} \rfloor ,
 2^n - 1 \}
 \end{equation}
 and for each such $m^{(1)}$ we have
 \begin{equation} \label{strata_2_bound}
 \lfloor \frac{ m^{(1)} }{2} \rfloor
 \leq m^{(12)}
 \leq min\{ \lfloor \frac{m^{(1)} + 2^{n-2} - 1}{2} \rfloor ,
 2^{n-1} - 1 \}.
 \end{equation}
Furthermore, all these bounds are attained.
\end{corollary}
\begin{remark} \label{remarkbound1flipnew}
 The lower bound for $m^{(1)}$ is precisely the inequality $ m -
 m^{(1)} \leq m^{(1)} + 1 $. Similarly,
 $ \lfloor \frac{ m^{(1)} }{2} \rfloor \leq m^{(12)} $ is
 equivalent to $ m^{(1)} - m^{(12)} \leq m^{(12)} + 1$.
\end{remark}
\begin{theorem} \label{simplecount}
 Let $m \in \Z$ with $2^n + 2^{n-1} -1 \leq m < 2^{n+1}$. Then
 $$
 a_m
 = \sum \{ a_{m^{(1)} }\, : \,
 m^{(1)} \, \textup{satisfies inequality~(\ref{strata_1_bound})} \}
 = \sum_{m^{(1)} = \lfloor \frac{m}{2} \rfloor}^{2^n -1}{ a_{m^{(1)}} }.
 $$
\end{theorem}
\begin{proof}
 Let $\lambda_0 + \lambda_1 + \cdots + \lambda_{n-1}$ be an
 $M$-partition of any such $m^{(1)}$.
 Since $2^n + 2^{n-1} - 1 \leq m$ then 
 $ 2^{n-1} \leq  m - (2^n - 1)  \leq   m - m^{(1)}$.
 By Lemma~\ref{bound1},
 $\lambda_{n-1} \leq \lceil \frac { m^{(1)} } {2} \rceil \leq
 2^{n-1} \leq m - m^{(1)}$. By Remark~\ref{remarkbound1flipnew},
 $ m - m^{(1)} \leq m^{(1)} + 1 $. Therefore, all partitions of
 $m^{(1)}$ satisfying inequality~(\ref{strata_1_bound}) extend to
 an $M$-partition of $m$ in the sense of Remark~\ref{extendedMremark}.
\end{proof}
\begin{example} \label{simplecountexample}
 The $M$-partitions of $25$ are extended from the $M$-partitions
 of $25^{(1)} = 12, 13, 14, 15$. Consequently,
 $a_{25} = a_{12} + a_{13} + a_{14} + a_{15}$. The $M$-partitions
 of $25$ are listed here with $25 - 25^{(1)}$ in bold.
\begin{align*}
 Mp(12) + {\bf 13} & = \{1+2+3+6 + {\bf 13}, \, 1+2+4+5 + {\bf 13} \}\\
 Mp(13) + {\bf 12} & = \{1+2+3+7 + {\bf 12}, \, 1+2+4+6 + {\bf 12}\}\\
 Mp(14) + {\bf 11} & = \{1+2+4+7 + {\bf 11}\}\\
 Mp(15) + {\bf 10} & = \{1+2+4+8 + {\bf 10} \}
\end{align*}
\end{example}

In general, for $2^n \leq m < 2^{n+1}$, not every $M$-partition of
an $m^{(1)}$ will have largest part no larger than $m - m^{(1)}$.
As a result, the calculation of $a_m$ may not be as
straightforward as that of Theorem~\ref{simplecount}.

\vspace{.2cm}

\noindent {\bf Example~\ref{countingexample} continued.}
 Let $m = 16$. By Corollary~\ref{bound1flipnew} we have
 $8 \leq 16^{(1)} \leq 11$. Thus the $M$-partitions are
 a subcollection of the following ordered compositions
 with $16 - 16^{(1)}$ in bold.
\begin{align*}
 Mp(8) + {\bf 8} \, = \, \{ & 1+1+2+4 +{\bf 8}, \, 1+1+3+3 +{\bf 8},
 \, 1+2+2+3 +{\bf 8} \}\\
 Mp(9) + {\bf 7} \, = \, \{ & 1+1+2+5+ {\bf 7}, \, 1+1+3+4+ {\bf 7},
 \, 1+2+2+4+ {\bf 7}, \\
                   & 1+2+3+3+ {\bf 7}\}\\
 Mp(10) + {\bf 6} \, = \,  \{ & 1+1+3+5 + {\bf 6}, \, 1+2+2+5+ {\bf 6},
 \, 1+2+3+4+ {\bf 6}\}\\
 Mp(11) + {\bf 5} \, = \, \{ & \underline{1+1+3+6+ {\bf 5}}, \,
 \underline{1+2+2+6+ {\bf 5}}, \, 1+2+3+5+ {\bf 5}, \\
                  & 1+2+4+4+ {\bf 5}\}
\end{align*}
 The two underlined compositions are not partitions because of the
 order on their parts but they do have the same parts as the
 compositions directly above them and these are $M$-partitions.
 Excluding the two underlined compositions, the remaining $12$
 ordered compositions are $M$-partitions and so $a_{16} = 12$.

\vspace{.2cm}

In the proof of Theorem~\ref{simplecount}, $2^n + 2^{n-1} -1 \leq
m$ was only required for $\lambda_{n-1} \leq m - m^{(1)}$. All the
other conditions of Remark~\ref{extendedMremark} were honored by
virtue of inequality~(\ref{strata_1_bound}). Keeping in mind that
the $M$-partitions of $m^{(1)}$ are constructed on $M$-partitions
of $m^{(12)}$ satisfying inequality~(\ref{strata_2_bound}), we can
once again re-interpret Remark~\ref{extendedMremark} as follows.
\begin{remark} \label{strataremark}
 The number of $M$-partitions of $m$, $a_m$ equals the cardinality of
 the set of partitions given by \newline
 \noindent
 $ \M_1 \, := \, \{ \,
  \lambda_0 + \lambda_1 + \cdots + \lambda_{n-2} + (m^{(1)} - m^{(12)} )
  \,:
  \, \lambda_{n-2} \leq m^{(1)} - m^{(12)} \leq m - m^{(1)}$ \newline
  \hspace{.4cm} 
  and
  $\lambda_0 + \lambda_1 + \cdots + \lambda_{n-2} +(m^{(1)} -
  m^{(12)})$ is an $M$-partition of $m^{(1)} \}$.
\end{remark}
Next we have a simple lemma that characterizes those partitions of
$m^{(1)}$ that do not extend to $M$-partitions of $m$.
\begin{lemma} \label{twoisenough}
 Let $m \in \Z$ with $2^n \leq m < 2^{n+1}$ and assume that
 $m^{(1)} - m^{(12)} > m - m^{(1)}$.
 If $m^{(12)} = \lambda_0 + \lambda_1 + \cdots + \lambda_{n-2}$
 is an $M$-partition then
 $\lambda_{n-2} < m^{(1)} - m^{(12)}$.
\end{lemma}
\begin{proof}
 Since $ m^{(12)} < 2^{n-1}$ then
 $3m^{(12)} = 2m^{(12)} + m^{(12)}
 \leq 2(2^{n-1}-1) + m^{(12)}
 = 2^n + m^{(12)} - 2$. Also, $2^n \leq m$
 which implies that $3m^{(12)} \leq  m + m^{(12)} -2 $.
 By assumption we have $  m + m^{(12)} < 2m^{(1)}$
 and so $3 m^{(12)} < 2 m^{(1)} - 2$.
 Subtracting $2 m^{(12)} - 2$ from both sides yields
 $ \frac{m^{(12)} + 2}{2} <  m^{(1)} - m^{(12)}$.
 But $ \lceil \frac{m^{(12)} }{2} \rceil < \frac{m^{(12)} +
 2}{2}$ and so
 $ \lceil \frac{m^{(12)} }{2} \rceil < m^{(1)} - m^{(12)}$.
 Since $\lambda_{n-2}$ is the largest part of an $M$-partition of
 $ m^{(12)}$ then, by Lemma~\ref{bound1}, we have
 $\lambda_{n-2} \leq \lceil \frac{ m^{(12)} }{2} \rceil$
 and so $ \lambda_{n-2} < m^{(1)} - m^{(12)}$.
\end{proof}

We will now calculate the cardinality of the set $Mp(m)$ by
determining the cardinality of the set $\M_1$ described in
Remark~\ref{strataremark}. We will do so by a recurrence relation.
\begin{theorem} \label{recurrence}
 For any $m \in \Z$ there is a recurrence relation for
 $a_m$ given by
 $$
 a_m \, = \,
 \sum_{m^{(1)} \, = \, \lfloor \frac{m}{2} \rfloor}
     ^{ min \{ \lfloor \frac{m + 2^{n-1} - 1}{2} \rfloor , \, 2^n - 1 \} }
      { \{ \, a_{m^{(1)}} \, - \,
              {\sum_{m^{(12)} = \lfloor \frac{m^{(1)}}{2} \rfloor}
                   ^{2m^{(1)} - m - 1}{a_{m^{(12)}  } } }  \, \} }.
 $$
\end{theorem}
\begin{proof}
 Let $\M$ equal the set
 $ \{ \, \lambda_0 + \lambda_1 + \cdots +
 \lambda_{n-2} + (m^{(1)} - m^{(12)} )  \,: \,
 \lambda_0 + \lambda_1 + \cdots + \lambda_{n-2} +(m^{(1)} -
 m^{(12)})\,
 \textup{is an $M$-partition of} \, m^{(1)} \}$
 and let
 $\M_2$ equal the subset of $\M$ given by
 $  \{ \,
 \lambda_0 + \lambda_1 + \cdots + \lambda_{n-2} + (m^{(1)} - m^{(12)} )
 \,: \, m^{(1)} - m^{(12)} > m - m^{(1)}
 \, \textup{and} \,
 \lambda_0 + \lambda_1 + \cdots + \lambda_{n-2} +(m^{(1)} -
 m^{(12)})\, \textup{is an $M$-partition of} \, m^{(1)} \}
 $.
 Then $\M = \M_1 \cup \M_2$ is a disjoint union of the set $\M$ and so
 we must have $ | \M_1 | \, = \, | \M | \, - \, | \M_2 | $.

 The set $\M$ is the set of all $M$-partitions of $m^{(1)}$ satisfying
 inequality~(\ref{strata_1_bound}).
 $$
 \M \,= \,
 \bigcup ^ {min\{ \lfloor \frac{m + 2^{n-1} - 1}{2} \rfloor , 2^n - 1 \}}
         _{m^{(1)}= \lfloor \frac{m}{2} \rfloor}
         { \, Mp(m^{(1)}) } \, .
 $$
 On the other hand, Lemma~\ref{twoisenough} says that $\M_2$ is in bijection
 with the set of all $M$-partitions of $m^{(12)}$ for all $m^{(12)}$
 satisfying inequality~(\ref{strata_2_bound}) and with
 $m^{(1)} - m^{(12)} > m - m^{(1)}$. That is, $\M_2$ is in bijection with
 the set of all $M$-partitions of $m^{(12)}$ with
 $\lfloor \frac{m^{(1)}}{2} \rfloor \leq m^{(12)} \leq
 2m^{(1)} - m - 1$ where $m^{(1)}$ satisfies
 inequality~(\ref{strata_1_bound}) and thus we can write the cardinality of
 $\M_2$ as
 $$
  | \, \M_2 \, | \, = \,
  \big| \,
 \bigcup ^ {min\{ \lfloor \frac{m + 2^{n-1} - 1}{2} \rfloor , 2^n - 1 \}}
         _{m^{(1)}= \lfloor \frac{m}{2} \rfloor}
 { \bigcup _{m^{(12)}= \lfloor \frac{m^{(1)}}{2} \rfloor}
            ^{2m^{(1)}-m-1}
            {Mp(m^{(12)})} } \, \big| \, .
 $$
 Recalling that $a_m \, = \,  | \, \M \, | - | \, \M_2 \, |$, we can write
 $$
 a_m   \, = \,
 \sum_{m^{(1)} \, = \, \lfloor \frac{m}{2} \rfloor}
     ^{ min \{ \lfloor \frac{m + 2^{n-1} - 1}{2} \rfloor , \, 2^n - 1 \} }
      { \{ \, a_{m^{(1)}} \, - \,
              {\sum_{m^{(12)} = \lfloor \frac{m^{(1)}}{2} \rfloor}
                   ^{2m^{(1)} - m - 1}{a_{m^{(12)}  } } }  \, \} }.
 $$
\end{proof}
\begin{table} [b]
 \begin{center}
  \begin{tabular}{||l|l||l|l||l|l||l|l||} \hline
$m$     &$a_m$  &$m$    &$a_m$  &$m$    &$a_m$  &$m$    &$a_m$ \\
\hline
1       &1      &17     &15     &33     &91     &49     &26     \\
        \hline
2       &1      &18     &13     &34     &82     &50     &20     \\
        \hline
3       &1      &19     &14     &35     &89     &51     &20     \\
        \hline
4       &1      &20     &11     &36     &77     &52     &14     \\
        \hline
5       &2      &21     &12     &37     &80     &53     &14     \\
        \hline
6       &1      &22     &9      &38     &70     &54     &10     \\
        \hline
7       &1      &23     &10     &39     &73     &55     &10     \\
        \hline
8       &3      &24     &6      &40     &60     &56     &6      \\
        \hline
9       &4      &25     &6      &41     &63     &57     &6      \\
        \hline
10      &3      &26     &4      &42     &53     &58     &4      \\
        \hline
11      &4      &27     &4      &43     &54     &59     &4      \\
        \hline
12      &2      &28     &2      &44     &43     &60     &2      \\
        \hline
13      &2      &29     &2      &45     &44     &61     &2      \\
        \hline
14      &1      &30     &1      &46     &35     &62     &1      \\
        \hline
15      &1      &31     &1      &47     &36     &63     &1      \\
        \hline
16      &12     &32     &84     &48     &26     &64     &908    \\
        \hline

   \end{tabular}
   \vspace{.2cm}
  \caption{Values of $a_m$ for $ 1 \leq m \leq 64$.}
 \end{center}
\end{table}
\begin{remark}
 As we would expect, Theorem~\ref{simplecount} follows as a special case
 of Theorem~\ref{recurrence}. The set $\M_2$ contains no elements precisely
 when $\lfloor \frac{m^{(1)}}{2} \rfloor > 2 m^{(1)} - m - 1$. This
 occurs only if $2^n + 2^{n-1} -1 \leq m < 2^{n+1}$.
\end{remark}

\vspace{.2cm}

\noindent {\bf Example~\ref{countingexample} continued.}
 Let $m = 16$. According to Theorem~\ref{recurrence}
 $a_{16} = a_8 + a_9 + a_{10} + (a_{11} - a_5) \, = \,
           3 + 4 + 3 + (4-2) $.
 The only instance of
 $\lfloor \frac{m^{(1)}}{2} \rfloor \leq 2 m^{(1)} - m - 1$
 being satisfied is when $16^{(1)} = 11$. Looking at the
 $M$-partitions of $11$ we see that there are two of them with
 largest part larger than $16 -11 = 5$; $\underline{1+1+3+6}$
 and $ \underline{1+2+2+6}$. Both of these $M$-partitions
 have largest part $6$ and so must be built upon all the
 $M$-partitions of $5$. Hence, we subtract $a_5$ from $a_{11}$.

\section{Simplifying the recurrence relation for $|Mp(m)|$}
\label{simplifyingsection}

In this section we exhibit a generating function for $m \in \Z$
provided that $2^n + 2^{n-1} - 1 \leq m  < 2^{n+1}$. In particular
the recurrence relation of Theorem~\ref{simplecount} has a
generating function.
\begin{lemma} \label{even}
 For even $m$ with $2^n + 2^{n-1} \leq m  < 2^{n+1} $ we have
 $a_m = a_{m+1}$.
\end{lemma}
\begin{proof}
 Since $ 2^n + 2^{n-1} \leq m \leq 2^{n+1} - 2$ then
 Theorem~\ref{simplecount} will suffice to calculate both $a_m$
 and $a_{m+1}$. Since $m$ is even then
 $\lfloor \frac{m}{2} \rfloor = \lfloor \frac{m+1}{2} \rfloor$
 and so the recurrence relation of Theorem~\ref{simplecount} is
 the same for both $a_m$ and $a_{m+1}$.
\end{proof}

We prove another lemma which will play a crucial role in the proof
of the main theorem of this section.
\begin{lemma} \label{brecurrence}
 For any integer $j \geq 0$ define the recurrence relation
 $ b_j = b_{j-1} + b_{\lfloor \frac{j}{2} \rfloor}$ with initial
 condition $b_0 = 1$.
 Then
 $
 b_j = \sum_{i=0}^{j}{b_{\lfloor \frac{i}{2} \rfloor}}
 $
\end{lemma}
\begin{proof}
 The lemma is true for $j=0$ and $j=1$. Utilizing an induction argument
 assume true for all $j < l$. Then
 $$
 b_l =b_{l-1} + b_{\lfloor \frac{l}{2} \rfloor} =
 \sum_{i=0}^{l-1}{b_{\lfloor \frac{i}{2} \rfloor}} +
 b_{\lfloor \frac{l} {2} \rfloor}
 =  \sum_{i=0}^{l}{b_{\lfloor \frac{i}{2} \rfloor}}.
 $$
 The last two equalities follow by the inductive hypothesis and so
 our claim is true for any non negative integer.
\end{proof}

The recurrence relation of Lemma~\ref{brecurrence} provides a more
efficient accounting of $a_m$ than that in
Theorem~\ref{simplecount}.
\begin{lemma}
 Let $m \in \Z$ satisfying
 $ 2^n + 2^{n-1} - 1 \leq m \leq 2^{n+1} - 1$
 and write $m$ in the form $m = 2^{n+1} - 1 -k$. Then
 $a_m = a_{2^{n+1} - 1 -k} = b_ {\lfloor \frac{k}{2} \rfloor}$.
\end{lemma}
\begin{proof}
 This will be shown by induction on $n = \lfloor \log_2 m \rfloor$.
 From Table~1 we can see that our claim is true for
 $ n = 0,1,2$ and so assume that our claim is true for all positive
 integers less than some $n$ and pick an $m$ such that
 $ 2^n + 2^{n-1} - 1 \leq m \leq 2^{n+1} - 1$.

 Since $m = 2^{n+1} - 1 -k$ then
 $\lfloor \frac{m}{2} \rfloor
 = \lfloor \frac {2^{n+1} - 2  - (k-1)}{2} \rfloor
 = 2^n - 1  - \lfloor \frac{k}{2} \rfloor$.
 From Theorem~\ref{simplecount} we have
 $$
 a_m \, =
 \sum_{m^{(1)}= \lfloor \frac{m}{2} \rfloor}
     ^{2^n - 1}
     {a_{m^{(1)}}}
 \, =
 \sum_{m^{(1)} = 2^n - 1  - \lfloor \frac{k}{2} \rfloor}
     ^{2^n - 1}
     {a_{m^{(1)}}}.
 $$
 Each $m^{(1)}$ satisfies $\lfloor \log_2 m^{(1)} \rfloor = n-1$
 and so our inductive hypothesis says this last summand (after reversing
 the order of summation) can be expressed as follows
 $$
 a_m \,
 = \sum_{ i  = 0 } ^ { \lfloor \frac{k}{2} \rfloor}
   { a_{2^{n} - 1 - i }} \, \,
 = \sum_{i=0}^{\lfloor \frac{k}{2} \rfloor}
        {b_{\lfloor \frac{i}{2} \rfloor}}
 \, = \, b_{\lfloor \frac{k}{2} \rfloor}.
 $$
 The last equality comes from Lemma~\ref{brecurrence}.
\end{proof}

We say a partition is {\em binary} if all its parts are powers of
$2$. See~\cite{partitions} for results about such partitions.
In~\cite{knuth}, Knuth studied binary partitions whose parts were
all distinct and, amongst other things, derived the following
result.
\begin{theorem} \label{knuththeorem}(Knuth)
 The number of distinct binary partitions of $2j$ into powers of
 $2$ equals $b_j$ where $b_j$ is the recurrence relation in
 Lemma~\ref{brecurrence}. Furthermore, this recurrence relation has
 a generating function given by
 $$
 (1-x)^{-1}\prod_{j=0}^{\infty}{(1-x^{2^{j}})^{-1}}.
 $$
\end{theorem}
We consequently have our main result which is a generating
function for $a_m$ when $2^n + 2^{n-1} - 1 \leq m < 2^{n+1}$.
\begin{corollary} \label{genfunction}
 If  $2^n + 2^{n-1} - 1 \leq m \leq 2^{n+1} - 1$
 and $m = 2^{n+1} - 1 - k$
 then $a_m$ equals the coefficient of
 $x^{ \lfloor \frac{k}{2} \rfloor }$
 in the above generating function.
\end{corollary}
 For the case of $2^n \leq m \leq 2^n + 2^{n-1} - 2$ it appears
 that the best we can do is the following:
 If $2^n \leq m = 2^{n+1} - 1 - k \leq 2^n + 2^{n-1} - 2$
 and $m^{\prime} = 2^n - 1 - k^{\prime}$ is such that
 $m = 2^{n-1} + 2^{n-2} + m^{\prime}$ then
 $ b_{ \lfloor \frac{k}{2} \rfloor} - a_m
 = b_{ \lfloor \frac{k^{\prime}}{2} \rfloor} - a_{m^{\prime}}$.
 However, it seems that no generating function can be arrived at
 for $m$ in the interval $2^n \leq m \leq 2^n + 2^{n-1} - 2$. In
 other words, a generating function for the recurrence relation
 of Theorem~\ref{recurrence} could not be arrived at.


\section*{Acknowledgements}

 I wish to express my deepest gratitude to Professors Hendrik W.
 Lenstra, Jr. and Bjorn Poonen for their encouragement, advice and
 direction during the research of this paper. 
 I thank George Andrews for referring me to the work of MacMahon. 
 Katia Hayati helped with the calculation for Table~1 and 
 John Sullivan was a source of good suggestions. 
 Rekha Thomas, Sara Billey and Peter Couperus provided valuable feedback 
 in the closing stages of writing this paper.


\begin{thebibliography}{2}

\bibitem{partitions} {\sc G. E. Andrews} , The Theory of
 Partitions, {\em Addison-Wesley} (1976)

\bibitem{knuth} {\sc D.E. Knuth}, An almost linear recurrence
 {\em Fibonacci Quarterly} {\bf 4}, (1966), 117--128.

\bibitem{macmahon} {\sc P. A. MacMahon}, The theory of perfect
 partitions and the compositions of multipartite numbers,
 {\em Messenger of Math.} {\bf 20}, (1891), 103--119

\bibitem{sloane} {\sc S. Plouffe and N.J.A. Sloane},
 The Encyclopedia of Integer Sequences {\em Academic Press} (1995)

\end{thebibliography}
\end{document}